%
%
     \input amstex
     \documentstyle{amsppt}
     \nopagenumbers

\pageheight{9 in}
\pagewidth{6.5 in}

\def\today{\ifcase\month\or January\or February\or March\or April\or
May\or June\or July\or August\or September\or October\or
November\or December\fi \space\number\day, \number\year }
%
     \NoBlackBoxes
     \nologo
     \parskip=3 pt
     \parindent = 0.3 true in
     \TagsOnRight
\def\ZZ{\Bbb Z}
\def\QQ{\Bbb Q}
\def\Q{\Bbb Q}

\def\epsilon{\varepsilon}
\def\phi{\varphi}
\def\pp{\frak p}

\def\mod{\,\text{mod}\,}

     \def\cf{{cf.}\ }
     \def\eg{{e.g.,}\ }

\def\EK4{E_{\! K,4} }

%
\count1=0
\count2=0

\document

    \footnote""{last modified \today}  

\footnote""{2010 {\it Mathematics Subject Classification.} Primary 15B35, 11R11, 11Y55 ; Secondary 05B20}
\footnote""{{\it Keywords}: Quadratic residue, sign matrix, power residue, reciprocity laws, splitting configurations.}

\topmatter

\title 
Characterizations of quadratic, cubic, and quartic residue matrices
\endtitle

\author
David S\. Dummit, Evan P\. Dummit, Hershy Kisilevsky
\endauthor

\abstract
We construct a collection of matrices defined by quadratic residue symbols, termed ``quadratic residue matrices'', 
associated to the splitting behavior of prime ideals in
a composite of quadratic extensions of $\QQ$, and prove 
a simple criterion characterizing such matrices. 
We also study the analogous classes of matrices constructed from the cubic and quartic residue symbols 
for a set of prime ideals of $\QQ(\sqrt{-3})$ and $\QQ(i)$, respectively. 
\endabstract

\address
University of Vermont, Dept. of Mathematics and Statistics, 16 Colchester
Ave., Burlington, VT 05401 
\endaddress

\address
University of Rochester, Dept. of Mathematics, 915 Hylan Building,
RC Box 270138, Rochester, NY 14627 
\endaddress

\address
Concordia University, Dept. of Mathematics and Statistics, 1400 De
Maisonneuve W., Montreal, QC
\endaddress

\email
dummit@math.uvm.edu
\endemail

\email
edummit@ur.rochester.edu
\endemail

\email
hershy.kisilevsky@concordia.ca
\endemail



\endtopmatter

\noindent
{\bf 1. Introduction}

\medskip
Let $n>1$ be an integer and $p_{1}$, $p_{2}$, $\dots$, $p_{n}$
be a set of distinct odd primes. 

The possible splitting behavior of the $p_{i}$ in the
composite of the quadratic extensions $\Q(\sqrt{p_{j}^*})$, where
$p^* = (-1)^{(p-1)/2} p$ (a minimally tamely ramified multiquadratic
extension, \cf [K-S]), is determined by the quadratic residue symbols 
$\fracwithdelims(){p_i}{p_j}$.  Quadratic reciprocity imposes
a relation on the splitting of $p_{i}$ in $\Q(\sqrt{p_{j}^*})$ and
$p_{j}$ in $\Q(\sqrt{p_{i}^*})$ and this leads to the definition of 
a ``quadratic residue matrix''.  The main purpose of this article is
to give a simple criterion that characterizes such matrices.  These matrices
seem to be natural elementary objects for study, 
but to the authors' knowledge have not previously appeared in the literature. 

We then consider higher-degree variants of this question arising from cubic and quartic 
residue symbols for primes of $\QQ(\sqrt{-3})$ and $\QQ(i)$, respectively.

\bigskip
\noindent
{\bf 2. Quadratic Residue Matrices}

\medskip
We begin with two elementary definitions.

\bigskip
\noindent
{\bf Definition.}
A ``sign matrix'' is an $n \times n$ matrix whose diagonal entries are all 0
and whose off-diagonal entries are all $\pm 1$.

\medskip
Among the sign matrices are those arising from Legendre symbols for a collection of
odd primes:

\medskip
\noindent
{\bf Definition.}
The ``quadratic residue'' (or ``QR'') matrix associated to the odd primes $p_{1}$, $p_{2}$,
$\dots$, $p_{n}$ is the $n\times n$ matrix whose $(i,j)$-entry
is the quadratic residue symbol $\fracwithdelims(){p_{i}}{p_{j}}$.

For example, for the primes $p_{1}=3$, $p_{2}=7$, and $p_{3}=13$,
the associated quadratic residue matrix is 
$$
M=\pmatrix
\format \r & \quad \r & \quad \r \\
0 & -1 & 1\\
1 & 0 & -1\\
1 & -1 & 0
\endpmatrix
$$

Our goal is to characterize those sign matrices that are QR matrices, as it is not immediately 
obvious whether, for example, the matrix 
$
M=\pmatrix
\format \r & \quad \r & \quad \r \\
0 & -1 & -1\\
-1 & 0 & -1\\
1 & 1 & 0
\endpmatrix
$
arises as the QR matrix for some set of primes.

Observe that the property of whether a matrix is a sign matrix or a quadratic residue matrix is 
invariant under conjugation by a permutation matrix, since the diagonal entries and underlying primes
will each be permuted.  If $p_1, \dots, p_n$ are odd primes and we order
them so that the first $s$, $1 \le s \le n$, are primes congruent to 3 modulo 4 and the
remaining $n-s$ are congruent to 1 modulo 4, then by quadratic reciprocity 
the associated quadratic residue matrix $M$ has the form
$$
\pmatrix
A & B \\
B^t & S
\endpmatrix
$$
where $A$ is an $s \times s$ skew-symmetric sign matrix, $S$ is an $(n-s) \times (n-s)$ symmetric
sign matrix and $B$ is an $s\times(n-s)$ matrix all of whose entries are $\pm1$, and $B^t$ denotes
the transpose of $B$.  

The following result proves the converse of this holds and also provides a simple criterion to determine when a
given sign matrix arises as a quadratic residue matrix for some set of primes.

\medskip
\noindent
{\bf Theorem 1.} If $M$ is an $n \times n$ sign matrix, then the following are equivalent:

\item {(a)}
There exists an integer $s$ with $1 \le s \le n$ such that the matrix $M$ can be 
conjugated by a permutation matrix into a block matrix of the form
$$
\pmatrix
A & B \\
B^t & S
\endpmatrix
$$
where $A$ is an $s \times s$ skew-symmetric sign matrix, $S$ is an $(n-s) \times (n-s)$ symmetric
sign matrix and $B$ is an $s\times(n-s)$ matrix all of whose entries are $\pm1$.  Here $B^t$ denotes
the transpose of $B$.

\item {(b)}
The matrix $M$ is the quadratic residue matrix associated to
a set of odd primes $p_{1}, p_{2}, \dots, p_{n}$.

\item {(c)}
There exists an integer $s$ with $1 \leq s \leq n$ such that the diagonal entries of $M^{2}$
consist of $s$ 
occurrences
of $n+1-2s$ and $n-s$ 
occurrences
of $n-1$.

\medskip
\noindent
{\smc Proof:} 
(a) implies (b): Suppose that $M = (m_{i,j})$ is a block matrix as in (a).  
We inductively construct primes $p_{1}, \dots , p_{n}$ for which
$M$ is the quadratic residue matrix: for the base case, let $p_{1}$ be any
prime congruent to 3 modulo 4. For the inductive step, suppose that
$p_{1},$ ... , $p_{k}$ are primes such that $\left(\dfrac{p_{i}}{p_{j}}\right)=m_{i,j}$
for $1\leq i,j\leq k$. For each $1\leq j\leq k$, choose a nonzero
residue class $u_{j}$ modulo $p_{j}$ such that $\left(\dfrac{u_{j}}{p_{j}}\right)=m_{k+1,j}$.
By the Chinese Remainder Theorem and Dirichlet's Theorem on primes
in arithmetic progression we may  choose a prime $p_{k+1}$ satisfying the congruences
$p_{k+1}\equiv u_{j}$ (mod $p_{j}$), along with either the congruence
$p_{k+1}\equiv 3$ (mod 4) if $k+1\leq s$ or $p_{k+1}\equiv 1$ (mod
4) if $k+1>s$. By construction, we have $\left(\dfrac{p_{k+1}}{p_{j}}\right)=m_{k+1,j}$
for all $1\leq j\leq k$, and quadratic reciprocity along with the
form of $M$ ensures that also $\left(\dfrac{p_{i}}{p_{k+1}}\right)=m_{i,k+1}$
for all $1\leq i\leq k$ is satisfied. Thus, $M$ is the quadratic residue matrix
associated to $p_{1}, \dots , p_{n}$, as claimed.

(b) implies (c):  
Suppose that $M$ is the quadratic residue matrix associated
to the primes $p_{1}$, $p_{2}$, $\dots$, $p_{n}$, and let $s$
be the number of these primes congruent to 3 modulo 4.
By rearranging
the primes, assume that $p_{1}$, $\dots$ $p_{s}$ are congruent
to 3 modulo 4 and that the remaining primes are congruent to 1 modulo
4. For $1\leq i\leq s$, the $i$th diagonal element of $M^{2}$ is
$$
(M^{2})_{i,i} = \sum_{j=1}^{n}\left(\dfrac{p_{i}}{p_{j}}\right)\left(\dfrac{p_{j}}{p_{i}}\right)=n+1-2s
$$
since by quadratic reciprocity the first $s$ terms are $-1$ (except
for the $i$th, which is 0), and the other $n-s$ terms are $+1$.
For $s+1\leq i\leq n$, the $i$th diagonal element of $M^{2}$ is
$$
(M^{2})_{i,i}=  \sum_{j=1}^{n}\left(\dfrac{p_{i}}{p_{j}}\right)\left(\dfrac{p_{j}}{p_{i}}\right)=n-1
$$
since by quadratic reciprocity all terms are $+1$ (again, except for the
$i$th, which is 0), proving (c).

(c) implies (a): 
Suppose that $M = (m_{i,j})$ is a sign matrix and 
the diagonal elements of $M^{2}$ consist of $s$ 
occurrences
of $n+1-2s$ and $n-s$ 
occurrences 
of $n-1$. By conjugating $M$ by an
appropriate permutation matrix, we may place the $s$ 
occurrences
of $n+1-2s$
in the first $s$ rows of $M^2$. For $s<i\leq n$, we have
$$
(M^{2})_{i,i}= \sum_{j=1}^{n}m_{i,j}m_{j,i}=n-1 ,
$$
but since there are only $n-1$ nonzero terms in the sum, we see that
$m_{i,j}=m_{j,i}$ for all $1\leq j\leq n$ and $s<i\leq n$. Then
for $1\leq i\leq s$, we have
$$
(M^{2})_{i,i}= \sum_{j=1}^{n}m_{i,j}m_{j,i}=n+1-2\cdot\#\{1\leq j\leq s\,:\, m_{i,j}=-m_{j,i}\}
$$
since $m_{i,j}m_{j,i}=+1$ whenever $j>s$, where we have also accounted
for the fact that $m_{i,i}=0$ in the count above. But now, since there are
at most $s$ terms in the count, and $(M^{2})_{i,i}=n+1-2s$, we see
that $m_{i,j}=-m_{j,i}$ for $1\leq j \leq s$.  Putting all of this together, we conclude that the matrix $M$ has
the block form in (a), completing the proof.

\medskip
\noindent
{\bf Example.}  The diagonal entries of the square of the QR matrix 
associated to the primes 3, 7, and 13 given previously are 0, 0, and 2, satisfying condition (c) of Theorem 1 with
$s = 2$.   On the other hand, the diagonal entries of the square of the sign matrix 
$
M=\pmatrix
\format \r & \quad \r & \quad \r \\
0 & -1 & -1\\
-1 & 0 & -1\\
1 & 1 & 0
\endpmatrix
$
mentioned earlier are 0, 0, and $-2$, showing this matrix is not a quadratic residue matrix.

\medskip
\noindent
{\bf Remark:}
It follows from (c) of Theorem 1 that $M$ is a QR matrix if and only if its transpose, $M^t$, is a QR matrix (and similarly for
$-M$ and $-M^t$), which is not immediately apparent from the definition.  

\medskip
\noindent
{\bf Remark:}
It also follows from Theorem 1 that we do not get a larger class of matrices than QR matrices 
by taking the elements of an $n \times n$ matrix to be the Jacobi symbol $\fracwithdelims(){P_{i}}{P_{j}}$ 
for a collection $P_1, \dots , P_n$ of odd, positive, pairwise relatively prime, integers (not necessarily prime): 
ordering the $P_i$ so the first $s$ are congruent to 3 modulo 4 and the remaining congruent to 
1 mod 4 and using quadratic reciprocity for the Jacobi symbol shows the matrix is in the form in 
(a) in the Theorem.

\bigskip
\noindent
{\bf 3. Counts for Quadratic Residue Matrices}

\bigskip

Computing the number and permutation equivalence classes of $n\times n$ QR matrices
by checking the criterion in (c) of Theorem 1 on all of the $n\times n$ sign matrices
rapidly becomes impractical (e.g., for $n > 5$) since the number of such sign matrices
is $2^{n(n-1)}$ and there are $n!$ permutation matrices.     

One can extend the computations slightly further (e.g., for $n = 6,7$) by using the
characterization in (b) of Theorem 1 and first computing the equivalence classes of 
symmetric and of skew-symmetric $m \times m$ sign matrices for $m \le n$.  This also rapidly becomes
impractical as there are $2^{m(m-1)/2}$ symmetric $m \times m$ sign matrices and the
same number of skew-symmetric $m \times m$ sign matrices.  The situation is further 
complicated by the fact that the ``reduced'' form $M = \pmatrix A & B \\ B^t & S \endpmatrix$
in the Theorem is not unique even after fixing the permutation equivalence classes
of $A$ and $S$.  This is already apparent when $n = 3$: there are 4 permutation equivalence
classes of $3 \times 3$ symmetric matrices (already in ``reduced'' form); such matrices are also in ``reduced'' form 
with $A = 0$ (the $1 \times 1$ skew-symmetric matrix) and $S$ one of the two permutation
inequivalent $2 \times 2$ symmetric
matrices, and there are 8 matrices of this type (since there are 4 possible choices for
the $1 \times 2$ matrix $B$)---in fact the 4 permutation equivalence classes of these 
8 matrices have sizes 1,1,3,3.  It is an interesting problem to determine whether there is 
a unique reduced form for QR matrices that would more readily allow the computation of the number of their
permutation equivalence classes.   

In addition to their use in the computation of QR matrices, the computation of the 
permutation equivalence classes of symmetric and of skew-symmetric
$n \times n$ sign matrices is a problem of independent interest.

The results of our numerical computations for the number of classes for sign
matrices are the following:

$$
\matrix
\underline{n}  & \underline{\text{no. of symmetric classes}} & \underline{\text{no. of skew-symmetric classes} } 
& \underline{ \text{number of matrices ($= 2^{n(n-1)/2}$)} } \\
2 & 2 & 1 & 2\\
3 & 4 & 2 & 8 \\
4 & 11 & 4 & 64 \\
5 & 34 & 12 & 1024 \\
6 & 156 & 56 & 32768 \\
7 & 1044 & 456 & 2097152
\endmatrix
$$
and for QR matrices the following:
$$
\matrix
\underline{n}  & \underline{\text{no. of QR matrix classes} } &  \underline{\text{no. of QR matrices}} &  
\underline{\text{no. of sign matrices ($= 2^{n(n-1)}$)}} \\
2 & 3 & 4 & 4\\
3 & 10 & 40 & 64 \\
4 & 47 & 768 & 4096 \\
5 & 314 & 27648  & 1048576 \\
6 & 3360 & 1900544 &  1073741824 \\
7 &  59744 & 253755392 & 4398046511104
\endmatrix
$$

There is also a graph-theoretic formulation of this counting question:
the number of $n \times n$ QR matrices is equal to the number of partially-directed
graphs on $n$ labeled vertices, such that (i) each vertex is colored
red or blue, (ii) between any two red points there is a single directed edge,
and (iii) between any red and blue or two blue points there is a single
undirected edge labelled with either {}``$+1$'' or {}``$-1$''.

The bijection between such graphs and QR matrices is as follows: red points
signify primes congruent to 3 modulo 4, and blue points signify 
primes congruent to 1 modulo 4, while an edge directed from
$p_{i}$ to $p_{j}$ indicates that 
$\fracwithdelims(){p_{i}}{p_{j}} = +1$
and 
$\fracwithdelims(){p_{j}}{p_{i}} = -1$, 
and the label on an undirected
edge joining $p_{i}$ and $p_{j}$ signifies the common value of
$\fracwithdelims(){p_{i}}{p_{j}}$ and $\fracwithdelims(){p_{j}}{p_{i}}$.

The sequence of numbers of permutation equivalence classes of symmetric sign matrices is sequence A000088 
in Sloane's database [S], defined there as the number of graphs on $n$ unlabeled nodes
(which is also the 
number of equivalence classes of sign patterns of totally nonzero symmetric $n \times n$ matrices), while the sequence 
of numbers of sign patterns of skew-symmetric sign matrices is sequence A000568, defined there as the number 
of outcomes of unlabeled $n$-team round-robin tournaments.  These interpretations follow immediately 
from the graph-theoretic version of the QR-matrix problem given above.
Explicitly, if the quadratic residue matrix is symmetric, all edges are labeled with $+1$ or $-1$;  the 
bijection with graphs on unlabeled nodes is given by deleting each edge labeled $-1$ and retaining each edge 
labeled $+1$.  If the quadratic residue matrix is skew-symmetric then all edges are directed;  the bijection 
with round-robin tournaments is given by viewing each directed edge as pointing from the loser of each match to the winner. 

The sequence of numbers of permutation equivalence classes of QR matrices and the sequence of numbers of
QR matrices do not seem to have been previously discovered.

\bigskip
\noindent
{\bf Remark:}
As mentioned in the Introduction, the authors were originally led to consideration of $n \times n$ 
QR matrices by a consideration of the possible decomposition types for the ramified primes in minimally tamely
ramified multiquadratic extensions of $\QQ$ with Galois group $(\ZZ/2\ZZ)^n$.  
There is a bijection between the possible decomposition
types of the ramified primes in such extensions and the QR matrices, and a bijection between the
various possible decomposition configuration ``types'' and the permutation equivalence classes of
QR matrices (\cf [D-K]).  

The computations here show, in particular, that there are 10 distinct types
of splitting of the three odd prime ideals $(p),(q),(r)$ in the extension 
$K = \QQ(\sqrt{p^*}, \sqrt{q^*},\sqrt{r^*} )$---for example where one of the prime ideals splits into precisely 4 prime ideals 
in $K$ and the other two prime ideals each splits into precisely 2 prime ideals in $K$.
Using the determination of the QR matrices, i.e., of the possible configuration types, one can then compute a frequency with which
each configuration occurs.  When $n = 3$, the 10 configuration types occur with frequencies
$\{1/32, 1/16, 1/16, 3/32, 3/32, 3/32, 3/32, 3/32, 3/16, 3/16 \}$.  For example, fields $K = \QQ(\sqrt{p^*}, \sqrt{q^*},\sqrt{r^*})$
where $(p),(q),(r)$ each split into 4 prime ideals in $K$ occur 1/32 of the time.

The counts of these decomposition types by direct computation (prior to our introduction of
QR matrices) did not immediately suggest any apparent frequencies of occurrence. 
For example, computing all 
306386 examples with  $p q r < 2457615 $ (so that each prime is less than the 15000th prime) led to frequencies 
$$
 \{0.037, 0.043, 0.062, 0.090, 0.108, 0.108, 0.123, 0.127, 0.138, 0.163 \}
$$
which, in spite of the relatively large sample size, do not compare particularly favorably with 
(to say nothing of actually suggesting) the correct frequencies of
$$
 \{0.03125, 0.0625, 0.0625, 0.09375, 0.09375, 0.09375, 0.09375, 0.09375 ,0.1875, 0.1875  \} .
$$
In fact, motivated by this example and as shown in [DGK], 
this discrepancy with the correct frequencies is not accidental, but rather due to a phenomenon
of large biases from small primes (similar to the ``prime races'' results comparing counts, for example, of primes
congruent to 1 and to 3 mod 4).  An important consequence of this result, reflecting the numerical
computations just mentioned, is that these biases are so strong that accurate computations of frequencies in 
problems of this type can be proved to be impossible---the computations required 
lie well beyond the limits of practicality. 
 
\bigskip
\noindent
{\bf 4. Cubic and Quartic Residue Matrices }

\bigskip
It is natural to consider generalizations of quadratic residue matrices to
matrices constructed from $m$th power residue symbols, in which case the
base field should contain the $m$th roots of unity.

\bigskip
\noindent
{\bf Definition.}
A ``cyclotomic sign matrix of $m$th roots of unity'' is an $n \times n$ matrix whose diagonal entries are all 0
and whose off-diagonal entries are all $m$th roots of unity.

\medskip
We consider here the cases $m = 3$ and $m = 4$ of cubic and quartic residues, and refer to the 
corresponding cyclotomic sign matrices simply as ``cubic sign matrices'' or ``quartic sign matrices'',
as appropriate. 
For larger values of $m$, the situation becomes more complicated, both because $m$th power reciprocity
becomes more complicated, and because prime ideals of the base field $\QQ(\zeta_m)$ of $m$th roots of unity 
need no longer be principal.  

\bigskip
\noindent


{\bf The case $\bold{m  = 3}$: Cubic Residue Matrices.}

In the case $m = 3$ of cubic residues, the most natural 
situation concerns the base field $\QQ( \sqrt {-3})$.  
The question that originally motivated our consideration
of QR matrices, the splitting of prime ideals in quadratic extensions of $\QQ$, in this context 
becomes a question of the splitting of prime ideals $\pp_1, \dots , \pp_n$ not dividing 3
of $K = \QQ( \sqrt {-3})$ in composites of cyclic cubic extensions of $K$.  If
$\pp$ is a prime ideal of $K$ not dividing 3, then $\pp$ is principal, and there
is a unique generator $\pi$ for $\pp$ which satisfies $\pi \equiv 1$ mod 3 in $K$
(a ``3-primary'' generator).   The minimally ramified cyclic cubic extensions of $K$ are then the Kummer extensions
$K (\root 3 \of \pi )$, the unique cubic subfield of the ray class field of conductor
$(\pp)$---the fact that $\pi$ is 3-primary ensures the extension is unramified at 3.  
From this perspective, the matrices to consider are given by the following.

\medskip
\noindent
{\bf Definition.}
The ``cubic residue matrix'' associated to the distinct prime ideals $\pp_1 , \dots, \pp_n$  of $\QQ( \sqrt {-3})$ not 
dividing 3 is the $n \times n$ matrix $M$ whose $(i,j)$-entry is the cubic residue symbol 
$\fracwithdelims(){\pi_{i}}{\pi_{j}}_3$ where $\pi_k$  is the
unique 3-primary generator for $\pp_{k}$ for $1 \le k \le n$.

\medskip
\noindent
The cubic residue symbol $\fracwithdelims(){\pi_{i}}{\pi_{j}}_3$ used in this Definition is the unique 3rd root of unity with
$$
\left(\dfrac{\pi_{i}}{\pi_{j}}\right)_3 \equiv \pi_i^{ (N\pi_j - 1 )/3 } \mod (\pi_j) \tag {1}
$$
where $N\pi_j$ denotes the norm from $\QQ( \sqrt {-3})$ to $\QQ$ of $\pi_j$. Recall also that
the cubic residue symbol $\fracwithdelims(){\pi_{i}}{\pi_{j}}_3$ gives the action of the
Frobenius automorphism $\sigma_j$ for the prime ideal $\pp_j = (\pi_j)$ in the cyclic cubic extension
$K (\root 3 \of \pi_i )$ of $K$:   $\sigma_j (\root 3 \of { \pi_i } ) = \fracwithdelims(){\pi_{i}}{\pi_{j}}_3 \root 3 \of {\pi_i} $
(\cf [C-F] for basic properties of $m$th power residue symbols).

For the cubic residue matrices the analogue of Theorem 1 is simpler and given by the following theorem.

\medskip
\noindent
{\bf Theorem 2.} If $M$ is an $n \times n$ cubic sign matrix (a matrix with 0's along the diagonal
and third roots of unity off the diagonal), then the following are equivalent:

\item {(a)}
The matrix $M$ is symmetric.

\item {(b)}
The matrix $M$ is the cubic residue matrix associated to
distinct prime ideals $\pp_1, \dots , \pp_n$ not dividing 3
in $\QQ( \sqrt {-3})$.

\medskip
\noindent
{\smc Proof:} 
If $\pi_i$ is a 3-primary generator for the prime ideal $\pp_i$, $i = 1, \dots , n$, then by cubic reciprocity we have 
$$
\left(\dfrac{\pi_{i}}{\pi_{j}}\right)_3 = 
\left(\dfrac{\pi_{j}}{\pi_{i}}\right)_3 
$$
(\cf, \eg  Exercise 2.14 in [C-F], or [I-R], p. 114), showing the cubic residue matrix $M$ associated to $\pp_1, \dots , \pp_n$
is symmetric.

Conversely, to see that every symmetric matrix $M$ whose diagonal entries are 0 and whose off-diagonal entries are
third roots of unity arises as such a matrix, we inductively construct prime ideals $\pp_1, \dots , \pp_n$ for
which $M$ is the associated cubic residue matrix.  Let $\pp_1$ be any prime ideal not dividing 3.   
For $l \ge 2$, suppose $\pp_1, \pp_2, \dots , \pp_{l-1}$
are prime ideals of $K = \QQ(\sqrt{-3})$ whose cubic residue symbols give rise to the first $(l-1) \times (l-1)$ upper left entries of $M$.   
As previously noted, specifying the third root of unity 
$\zeta = \fracwithdelims(){\pi_{i}}{\pp_{l}}_3$ is equivalent 
to specifying an element $\sigma$ in the Galois group of $K ( \root 3 \of \pi_i)$ over $K$ (by 
$\sigma  \root 3 \of \pi_i  = \zeta  \root 3 \of \pi_i ) $.  Since the extensions $K ( \root 3 \of \pi_i)$, 
$i = 1,2, \dots, l-1$ are linearly disjoint, it follows by 
Chebotarev's density theorem (applied to the composite of these extensions)
that there is a prime ideal $\pp_{l}$ whose cubic residue symbols agree
with the first $l$ elements in the $l$th column of $M$.  Then the symmetry of $M$ and cubic reciprocity
show that the first $l \times l$ upper left entries of $M$ are $\fracwithdelims(){\pi_{i}}{\pi_{j}}_3$ for the 3-primary
generators of the prime ideals $\pp_1, \pp_2, \dots , \pp_{l-1}, \pp_l$, showing inductively that $M$ is 
a cubic residue matrix.

\bigskip
While the number of such cubic residue matrices is then $3^{n(n-1)/2}$, the number of permutation
equivalence classes remains an interesting question.

\bigskip
\noindent
{\bf  The case $\bold{m  = 4}$: Quartic Residue Matrices.}

In the case $m=4$ of quartic residues, the most natural situation concerns the base field $K = \QQ( i)$.
In this case each prime ideal $\pp$ not dividing 2 has a
unique generator $\pi \equiv 1$ mod $2(1 + i)$ (a ``2-primary'' generator).  An element
$a + b i$ is 2-primary if either $a \equiv 1$ and $b \equiv 0$ modulo 4 or $a \equiv 3$ and $b \equiv 2$
modulo 4 (i.e., $a + b i$ is either $1$ or $3 + 2 i$ modulo (4) in $\QQ(i)$).  

The minimally tamely ramified cyclic quartic extensions of $K$ are the Kummer extensions
$K (\root 4 \of \pi )$, the unique cyclic quartic subfield of the ray class field of conductor
$(\pp)$---similar to the cubic case, the fact that $\pi$ is 2-primary ensures the extension is unramified at 2.  
The corresponding matrices to consider in this case are the following.

\medskip
\noindent
{\bf Definition.}
The ``quartic residue matrix'' associated to the distinct prime ideals $\pp_1 , \dots, \pp_n$  of $\QQ( i )$ not 
dividing 2 is the $n \times n$ matrix $M$ whose $(j,k)$-entry is the quartic residue symbol 
$\fracwithdelims(){\pi_{j}}{\pi_{k}}_4$ where $\pi_l$  is the
unique 2-primary generator for $\pp_{l}$ for $1 \le l \le n$.

\medskip
\noindent
The quartic residue symbol 
$\fracwithdelims(){\pi_{j}}{\pi_{k}}_4$ used in this Definition is the unique 4th root of unity with
$$
\left(\dfrac{\pi_{j}}{\pi_{k}}\right)_4 \equiv \pi_j^{ (N\pi_k - 1 )/4 } \mod (\pi_k) \tag {2}
$$
where $N\pi_k$ denotes the norm from $\QQ(i)$ to $\QQ$ of $\pi_k$.  As with cubic residue symbols,
the quartic residue symbol $\fracwithdelims(){\pi_{j}}{\pi_{k}}_4$ gives the action of the
Frobenius automorphism $\sigma_k$ for the prime ideal $\pp_k = (\pi_k)$ in the cyclic quartic extension
$K (\root 4 \of {\pi_j} )$ of $K$:   $\sigma_k (\root 4 \of { \pi_j } ) = \fracwithdelims(){\pi_{j}}{\pi_{jk}}_4 \root 4 \of {\pi_j} $.

Quartic reciprocity (see, for example, [I-R], p\. 123) can be written in the form
$$
\left(\dfrac{\pi_{j}}{\pi_{k}}\right)_4  \overline {\left(\dfrac{\pi_{k}}{\pi_{j}}\right)_4 } = (-1)^{\frac{N\pi_j - 1 }{4}  \frac{N\pi_k - 1 }{4} } ,
$$
where the bar denotes complex conjugation.
We also note one of the supplementary laws of quartic reciprocity that follows easily from (2), namely 
$ \fracwithdelims(){-1}{\pp}_4 = (-1)^{(a-1)/2}$ if $a + b i$ is the 2-primary generator for
the prime ideal $\pp$ not dividing 2. One consequence of this is that the 2-primary generator for the odd
prime ideal $\pp$ of $K = \QQ(i)$ is $1$ modulo (4) if and only if $\pp$ splits in the field 
$K(\root 4 \of {-1} ) = \QQ(i, \sqrt 2)$ of 8th roots of unity (and the 2-primary generator is
$3 + 2 i$ modulo (4) if and only if $\pp$ is inert in $\QQ(i, \sqrt 2)$ ).

The characterization of quartic residue matrices in the following theorem 
is similar to that of the quadratic residue matrices in Theorem 1, with the product
of $M$ and and its complex conjugate $\overline M$ taking the place of $M^2$.  

\medskip
\noindent
{\bf Theorem 3.} If $M$ is an $n \times n$ quartic sign matrix (a matrix with 0's along the diagonal
and fourth roots of unity off the diagonal), then the following are equivalent:

\item {(a)}
There exists an integer $s$ with $1 \le s \le n$ such that the matrix $M$ can be 
conjugated by a permutation matrix into a block matrix of the form
$$
\pmatrix
A & B \\
B^t & S
\endpmatrix
$$
where $A$ is an $s \times s$ skew-symmetric quartic sign matrix, $S$ is an $(n-s) \times (n-s)$ symmetric
quartic sign matrix and $B$ is an $s\times(n-s)$ matrix all of whose entries are $\pm1$ or $\pm i$.  


\item {(b)}
The matrix $M$ is the quartic residue matrix associated to
a set of distinct prime ideals $\pp_{1}, \pp_{2}, \dots, \pp_{n}$ not dividing 2 in $\QQ(i)$.

\item {(c)}
If $M = (m_{j,k})$, then $m_{j,k} = \pm m_{k,j}$ for all $j,k$ with $1 \le j,k \le n$, and there exists
an integer $s$ with $1 \leq s \leq n$ such that the diagonal entries of $M \overline M$
consist of $s$ 
occurrences
of $n+1-2s$ and $n-s$ 
occurrences
of $n-1$.

\medskip
\noindent
{\smc Proof:} 
(a) implies (b): Suppose $M=\{m_{j,k}\}$ is a block matrix as
in (a). We inductively construct prime ideals $\pp_{1}, \dots, \pp_{n}$
for which $M$ is the quartic residue matrix. For the base case, let $\pp_{1}$ be any
prime ideal of $\QQ(i)$ that is inert in the extension $\QQ( i, \sqrt 2)$, so as previously noted,
has 2-primary generator congruent to $3+2i$ modulo (4).

For the inductive step,
suppose that $\pp_{1}, \dots , \pp_{l-1}$ are distinct prime ideals not dividing 2 with 2-primary
generators $\pi_{1}, \dots , \pi_{l-1}$ satisfying
$\fracwithdelims(){\pi_{j}}{\pi_{k}}_4 = m_{j,k}$ for $1\le j,k \leq l-1$. 
Applying the Chebotarev density theorem in the extension 
$K(\sqrt 2 , \root 4 \of { \pi_1 } , \dots ,  \root 4 \of { \pi_{l-1} } )$ of $K$, there exists a
prime ideal $\pp_l$ of $K$ not dividing 2 whose Frobenius automorphism in the extension
$K(\root 4 \of { \pi_j } ) /K$ maps $\root 4 \of { \pi_j }$ to $m_{j,l}  \root 4 \of { \pi_j } $,
$1 \le j \le l-1$, 
and whose Frobenius automorphism in the extension $K(\sqrt 2 ) /K$ is trivial if $l > s$ and
is the nontrivial automorphism if $l \le s$.

Put another way, the 2-primary generator $\pi_l$ of $\pp_l$ satisfies 
$\pi_{l} \equiv 3+2i$ mod (4) if $l \le s$ and $\pi_{l} \equiv 1$ mod (4) if $l>s$, and
$\fracwithdelims(){\pi_{j}}{\pi_{l}}_4 = m_{j,l}$ for all $j$, $1\leq j\leq l$ . 

Since $\pi \equiv 1$ mod (4)  implies $N \pi \equiv 1 $ mod 8 and 
$\pi \equiv  3+2i$ mod (4) implies $N \pi \equiv 5 $ mod 8, it follows from 
quartic reciprocity, along with the form of $M$,  that
$\fracwithdelims(){\pi_{l}}{\pi_{j}}_4 = m_{l,j}$ for all $j$, $1\leq j\leq l$,
completing the induction and proving $M$ is a quartic residue matrix.

(b) implies (c): Suppose that $M$ is the quartic residue matrix associated
to the distinct prime ideals $\pp_{1}$, $\pp_{2} , \dots ,  \pp_{n}$ not dividing 2.
The first part of the criterion in (c) follows immediately from quartic reciprocity.

For the second part, rearrange the prime ideals, if necessary, so that 
$\pp_{1}, \dots , \pp_{s}$  have 2-primary generators
$\pi$ that are congruent to $3+2i$ mod (4) (i.e., have $N \pi \equiv 5 $ mod 8) and
the remaining prime ideals have 2-primary generators
$\pi$ that are congruent to $1$ mod (4) (i.e., have $N \pi \equiv 1 $ mod 8).
For $1 \leq j \leq s$, the $j$th diagonal element of $M\overline{M}$ is
$$
(M\overline{M})_{j,j}= \sum_{k=1}^{n} \left( \dfrac{\pi_{j}}{\pi_{k}} \right)_{4} \overline{ \left( \dfrac{\pi_{k}}{\pi_{j}} \right)_{4}} = n+1-2s 
$$
since by quartic reciprocity the first $s$ terms are $-1$ (except
for the $j$th, which is 0), and the other $n-s$ terms are $+1$.
For $s+1 \leq j \leq n$, the $j$th diagonal element of $M\overline{M}$ is 
$$
(M\overline{M})_{j,j}= \sum_{k=1}^{n} \left( \dfrac{\pi_{j}}{\pi_{k}} \right)_{4} \overline{ \left( \dfrac{\pi_{k}}{\pi_{j}} \right)_{4}} = n-1 
$$
since by quartic reciprocity all terms are $+1$ (except for the $j$th, which is 0), proving (c).

(c) implies (a): Suppose that $m_{j,k}=\pm m_{k,j}$ for each pair
$(j,k)$, and that the diagonal entries of the matrix $M \overline{M}$
consist of $s$ occurrences of $n+1-2s$ and $n-s$ occurrences of
$n-1$. 

By the assumption that $m_{j,k}=\pm m_{k,j}$ and these are 4th roots of
unity, the only possibility is for $m_{j,k} \overline{m_{k,j}}$ to be either $+1$ or $-1$.

By conjugating $M$ by an appropriate permutation matrix we may place the $s$ occurrences of $n+1-2s$
in the first $s$ rows of $M \overline{M}$. 
For $s < j \leq n$, we have 
$$
(M\overline{M})_{j,j}= \sum_{k=1}^{n} m_{j,k} \overline{m_{k,j}} = n-1 ,
$$
but since there are only $n-1$ nonzero terms in the sum, we see that
$m_{j,k} \overline{m_{k,j}}=1$ and hence $m_{j,k}=m_{k,j}$ for all
$1\leq k \leq n$ and $s < j \leq n$. 

For $1 \leq j \leq s$, we have
$$
(M\overline{M})_{j,j}= \sum_{k=1}^{n} m_{j,k} \overline{m_{k,j}} = n+1-2\cdot\#\{1 \leq k \leq s\,:\, m_{j,k} \overline{m_{k,j}} = -1 \}
$$
since $m_{j,k} \overline{m_{k,j}} = +1$ whenever $j>s$ and $m_{j,k}\overline{m_{k,j}}$
can only be $1$ or $-1$. But now since there at most $s$ terms
in the count, and $(M\overline{M})_{j,j}=n+1-2s$, we see that $m_{j,k}=-m_{k,j}$
for $1\leq k \leq s$, and $M$ has the form in (a), completing the proof.

\bigskip
\noindent
{\bf Remark:}
In defining the cubic and quartic residue matrices we have used ``primary'' generators of ideals 
in part because of the connection of these matrices with splitting questions in minimally ramified extensions.  
While it may not be so immediately apparent, the quadratic residue matrices are also constructed
the same way:  in the quadratic case, the `2-primary' generator of the odd prime ideal $(p)$ is the
element $p^* = (-1)^{(p-1)/2} p$, which leads to matrices whose entries are the quadratic residue symbols
$\fracwithdelims(){p_{i}^*}{p_{j}}$.  Since $\fracwithdelims(){p_{i}^*}{p_{j}} = \fracwithdelims(){p_{j}}{p_{i}}$
by quadratic reciprocity, these are just the transposes of the more elementary QR matrices we defined, 
which, as mentioned previously, also gives the QR matrices.   

\bigskip
\noindent
{\bf Remark:}
If we abandon the connection to splitting of prime ideals, the restriction to primary elements $\pi$ in the 
construction of $m$th power residue matrices above can be removed: Take a collection $\pi_1, \pi_2, \dots , \pi_n$ of prime
elements in the cyclotomic field $\QQ(\zeta_m)$ of $m$th roots of unity that are pairwise relatively prime and prime to $m$ 
and consider the $n \times n$ matrix
whose $(i,j)$ entry is the $m$th power residue symbol $\fracwithdelims(){\pi_{i}}{\pi_{j}}_m$, the unique $m$th root of unity congruent to
$\pi_i^{(N \pi_j - 1)/m}$ modulo the prime ideal $(\pi_j)$.  
Even for $m = 2 , 3, 4$ this gives a larger class of residue matrices than considered above: for example, when $m = 3$, these
matrices need not be symmetric;  when $m = 2$ these matrices are given using the Legendre symbols $\fracwithdelims(){p_{i}}{\vert p_{j} \vert}$
where the $p_i$ need not be positive, and the matrix 
$
\pmatrix
0 & -1 & 1 \\
1 & 0 & 1 \\
1 & -1 & 0 \\
\endpmatrix
$
associated to $-3,-5,7$ is not a QR matrix by Theorem 1(c).  Similarly, for $m = 2$, using the Kronecker symbol defines a larger class
(for example, the matrix associated to
$-3,-5,7$ in this case is
$
\pmatrix
0 & 1 & 1 \\
-1 & 0 & 1 \\
1 & -1 & 0 \\
\endpmatrix
$
and is again not a QR matrix---in fact every $3 \times 3$ sign matrix arises using the Kronecker symbol).
These other classes of matrices seem less tractable to characterization.

  \bigskip
  \noindent
  {\bf 5. Acknowledgements }
  
 The authors would like to thank Dinesh Thakur, John Doyle, Douglas Haessig, Amanda Tucker, and Andrew Bridy for their useful feedback and comments on the content of the paper.

\end